\documentclass[11pt]{article}
\usepackage{amsfonts}
\usepackage{dsfont}
\usepackage{amsfonts}
\usepackage{bbm}
\usepackage{mathrsfs}
\usepackage{color}
\usepackage{mathrsfs}
\usepackage{amsmath}
\usepackage{enumitem}
\allowdisplaybreaks[4]
\usepackage{amsfonts}
\usepackage{amsthm}
\usepackage{amssymb, amsmath, cite}
\usepackage{bm}
\usepackage{amsthm}
\usepackage{verbatim}
\usepackage{graphicx}
\usepackage{color}
\usepackage{amsfonts}
\usepackage{dsfont}
\usepackage{amsfonts}
\usepackage{bbm}
\usepackage{mathrsfs}
\usepackage{color}
\usepackage{graphicx,amsmath,amsfonts,amssymb,amscd,amsthm}
\usepackage{hyperref}
\usepackage[numbers,sort&compress]{natbib}
\usepackage{enumitem}
\setenumerate[1]{itemsep=1pt,partopsep=0pt,parsep=\parskip,topsep=1pt}
\setitemize[1]{itemsep=1pt,partopsep=0pt,parsep=\parskip,topsep=1pt}

\usepackage{graphicx,amsmath,amsfonts,amssymb,amscd,amsthm}
\usepackage{hyperref}
\usepackage{array}
\usepackage{booktabs}
\usepackage{multirow}
\usepackage{makecell}
\usepackage{abstract}
\hypersetup{colorlinks,linkcolor=blue,citecolor=blue}
\newtheoremstyle{kai}
{3pt} {3pt} {} {} {\bfseries} {.} {.5em} {}
\def\EquationsBySection{\def\theequation
{\thesection.\arabic{equation}}
\@addtoreset{equation}{section}}
\def\qed{\hfill\vrule width4pt height8pt depth0pt}

\newcommand\old[1]{}

 \newcommand{\Bp}{\begin{proof}}
 \newcommand{\Ep}{\end{proof}}
\renewcommand{\theequation}{\thesection.\arabic{equation}}

 \newcommand{\R}{\mathbb{R}}

 \newcommand{\beq}{\begin{equation}}
\newcommand{\eeq}{\end{equation}}

\newcommand{\be}{\begin{equation} \label}
\newcommand{\ee}{\end{equation}}
\newcommand{\bea}{\begin{eqnarray}\label}
\newcommand{\eea}{\end{eqnarray}}
\newcommand{\bas}{\begin{eqnarray*}}
\newcommand{\eas}{\end{eqnarray*}}
\newcommand{\bit}{\begin{itemize}}
\newcommand{\eit}{\end{itemize}}

\newtheorem{theo}{Theorem}[section]
\newtheorem*{theoa}{Theorem A}
\newtheorem*{theob}{Theorem B}
\newtheorem{lem}[theo]{Lemma}

\newtheorem{rem}[theo]{Remark}
\newtheorem{prop}[theo]{Proposition}

\numberwithin{equation}{section}
\addcontentsline{toc}{section}{Acknowledgement}
\begin{document}
\title{\bf Exact blow-up profiles for the parabolic-elliptic Keller-Segel system in dimensions $N\ge 3$}
\author{{\sc Xueli Bai} \thanks {School of Mathematics and Statistics, Northwestern Polytechnical University, 
Xi'an 710072, China; E-mail: xlbai2015@nwpu.edu.cn}
 \and {\sc Maolin Zhou}\thanks{Chern Institute of Mathematics and LPMC, Nankai University, Tianjin 300071, China; Corresponding author.
E-mail: zhouml123@nankai.edu.cn}}
\date{}
\maketitle
\begin{abstract}

In this paper, we obtain the exact blow-up profiles of solutions of the Keller-Segel-Patlak system in the space with dimensions $N\ge 3$:
\begin{equation}
    	\left\{ \begin{array}{ll}
	u_t = \Delta u - \nabla \cdot (u\nabla v),
	& x\in\R^N, \ t\in (0,T), \\[1mm]
	0 = \Delta v+u,
	& x\in\R^N, \ t\in (0,T), \\[1mm]
	u(x,0)=u_0(x),
	& x\in\R^N,
    	\end{array} \right.
\end{equation}
which solves an open problem proposed by P. Souplet and M. Winkler in \cite{Soup-Win}. To establish this achievement, we develop the zero number argument for nonlinear equations with unbounded coefficients and construct a family of auxiliary backward self-similar solutions through nontrivial ODE analysis.
\textbf{AMS (2020) Subject Classification}: {92C17, 35B40, 35B44, 35K40.}
\end{abstract}

\section{Introduction}

$\ \ \ \ $In this paper, we consider the blow-up profiles of solutions of the following Keller-Segel-Patlak system
\be{00}
    	\left\{ \begin{array}{ll}
	u_t = \Delta u - \nabla \cdot (u\nabla v),
	& x\in\R^N, \ t\in (0,T), \\[1mm]
	0 = \Delta v+u,
	& x\in\R^N, \ t\in (0,T), \\[1mm]
	u(x,0)=u_0(x),
	& x\in\R^N,
    	\end{array} \right.
\ee
where $T>0$ is the maximal existence time. This system can be considered as a simplified model of the classical Keller-Segel model which was firstly established by Keller and Segel in 1970s in \cite{KS1} and has been studied during the last four decades (see \cite{Aru-Tya,BBTW,BOSTW,HillP,Horstmann,Joh-Win}).

It is well-known that all the solutions of \eqref{00} are globally bounded for $N=1$ (see e.g. \cite{Osa-Yag}). When $N=2$, an interesting $8\pi$ critical mass phenomenon has been well studied about the blow-up solutions (see, e.g. \cite{Blan-Cal-Car1,Blan-Cal-Car2,cpam08,Bla-Dol-Per,Senba-Suzu,Jager-Luc,Nagai95,cpam18,Mizo22}).

For the case $ N\ge 3$, problem \eqref{00} admits both globally bounded solutions and \textcolor{blue}{blow-up} solutions, but the critical mass does not exist anymore (see e.g. \cite{Herrero-M-V2,Senba05,WinJDE,WinMathAnn}).  A family of backward self-similar blow-up solutions with profile $u(x,T)\sim\frac{C}{|x|^2}$ as $|x|\rightarrow 0$ were constructed by M.A. Herrero, E. Medina, J.J.L. Vel\'{a}zquez \cite{Herrero-M-V} (also see T. Senba \cite{Senba05} for a further investigation). This implies that $|x|^{-2}$ should be the order of the profiles of the blow-up solutions when $N\geq 3$.  Later, under appropriate conditions
\bea{initial 1}
0\le u_0\in L^\infty(\mathbb{R}^N),~u_0{\rm ~is~radially~symmetric~and~nonincreasing~with~respect~to~}  |x|,
\eea
and
\bea{initial 1'}
u_0\in C^1(\mathbb{R}^N), ~ u_0(x)\int_{|y|\le |x|} u_0(y)dy\geq |x|^{n-1}|\nabla u_0(x)| \mbox{ holds for }x\in \mathbb{R}^N,
\eea
P. Souplet and M. Winkler \cite{Soup-Win} showed that it really holds for general blow-up solutions of \eqref{00}.
\begin{theoa}[Souplet-Winkler 2019, \cite{Soup-Win}]  Assume $N\geq 3$ and $u_0$ satisfies \eqref{initial 1} and \eqref{initial 1'}. Let $u$ be the solution  of \eqref{00} with initial data $u_0$. Further assume that the blow-up time $T<\infty$ and the blow-up set $B(u)$ is not $\R^N$.
Then there exists $C\ge c>0$ such that
\be{upper of U}
c\le |x|^2U(x)\le C,~0<|x|\le 1,
\ee
where
$
U(x):=\lim_{t\rightarrow T}u(x,t)$ exists for all $x\in \R^N-\{0\}$.
\end{theoa}

In \cite{Soup-Win}, an open problem is naturally addressed that whether
the limit $\lim_{x\rightarrow 0}|x|^2 U(x)$ exists. In this paper, we completely solve this problem by developing an innovative  method inspired by the work of H. Matano and F. Merle \cite{Mat-Mer}.


\begin{theo}\label{thm1}
Assume $N\geq 3$ and $u_0$ satisfies \eqref{initial 1} and
\bea{initial 2}
\lim_{|x|\rightarrow \infty}|x|^{-n+2}\int_{|y|\le |x|}u_0(y)dy\mbox{ exists}.
\eea
Let $u$ be the solution  of \eqref{00} with initial data $u_0$. Further assume that the blow-up time $T<\infty$ and the blow-up set $B(u)$ is not $\R^N$.  Then we have
\begin{equation}\label{U profile}\lim_{|x|\rightarrow0^+}U(x)|x|^2=\alpha\in [0,\infty).
\end{equation}
\end{theo}
\begin{rem}
Compared with the literature \cite{Soup-Win},  the assumpution \eqref{initial 1'} is replaced by the condition \eqref{initial 2} due to technical reasons. In particular, a common hypothesis $u_0\in L^1$ upon biological background  is a special case of the condition \eqref{initial 2}.
\end{rem}
To better illustrate the main innovation in proving Theorem \ref{thm1}, we first briefly explain   the method employed in the proof of Theorem A in \cite{Soup-Win}. A crucial step  is to transform the blow-up problem from $u$ to $w$, where $w(r,t):=r^{-n}\int_0^r s^{n-1}u(s,t)$ is the average of $u$ in balls and satisfies
\begin{equation}\label{equ_w}
w_t=w_{rr}+\frac{N+1}{r}w_{r}+(nw+r w_r)w,
\end{equation}
where the comparison principle is valid.
Moreover, it is easy to show that $u$ and $w$  not only blow up simultaneously, but also have similar behaviors at the blow-up time. Hence, instead of $u$, the authors focus on studying  $w$ in \cite{Soup-Win}. Borrowing the idea of \cite{Fri-McL}, the upper bound of $w$ is  obtained  by showing that $J(r,t):=w_r-\varepsilon r w^2\ge 0$ holds in a fixed ball for small $\varepsilon>0$, while the lower bound of $w$ is derived by another method from \cite{Soup}. 
Therefore, in Theorem A, the connection between the upper and lower bounds is unclear and thus the sharp estimate on the blow-up profiles  remains as an open problem mentioned above.
To conquer this open problem, in this paper, we introduce a new way to combine zero number argument and the decay rate of backward self-similar solutions of equation \eqref{equ_w}. More specifically,  similar to the definition of $U$, we denote
\begin{equation}\label{W def}
W(r):=\lim_{t\rightarrow T}w(r,t)
\end{equation}
 And condtion \eqref{initial 2} is equivalent to
\be{initial w_0}
	\lim_{r\rightarrow\infty}w_0(r)r^2\in [0,\infty].
\ee
Under condition \eqref{initial w_0}, we can count the number of intersection points between $w$ and some auxilliary self-similar solution $\phi$ to show that the blow-up profile of $w$ can not osillate bewteen $\frac{C}{|x|^2}$ and $\frac{c}{|x|^2}$ with $c<C$ such that $ \lim_{|x|\rightarrow0^+}W(x)|x|^2$ exists and so is $U$.

Next, let us introduce our key idea heuristically. Inspired by \cite{Mat-Mer}, if we considered a blow-up phenomenon induced by the concentration of a soliton (or a self-similar soliton), the blow-up profile is just the concentration of the tail of the solution at infinity (after rescaling). This implies that it is very important to understand the decay rate of the backward self-similar solutions of equation (\ref{equ_w}) at infinity. Then we have the following result on a class of special backward self-similar solutions of equation \eqref{equ_w}.
\begin{prop}\label{dense S}
There exists a subset $\mathbb{S}\subset\mathbb{R}_+$ which is dense in $\mathbb{R}_+$ such that for any $m\in \mathbb{S}$, we can always find $\phi$ and $\ell\geq 0$ such that
$$\phi_{\xi\xi}+(\frac{N+1}{\xi}-\frac{\xi}{2})\phi_\xi-\phi+\phi(\xi \phi_\xi+N\phi)=0, \mbox{ for } \xi\in(\ell,\infty)\mbox{ and }
\lim_{\xi\rightarrow\infty}\phi(\xi)\xi^2=m.$$
Moreover, we can choose $\ell$ such that one of the following cases must hold
\begin{itemize}
\item $\phi$ is unbounded: $\phi(\xi)\rightarrow\infty$ as $\xi\rightarrow \ell^+$; \\
\item $\phi$ is bounded and $\ell>0$: $\lim\limits_{\xi\rightarrow \ell^+}\phi(r)=0$;\\
\item $\phi$ is bounded and $\ell=0$: $\lim\limits_{\xi\rightarrow \ell^+}\phi(\xi)$ exists and $\phi^\prime(\xi)\rightarrow 0$ as $\xi\rightarrow \ell^+$.
\end{itemize}
For convenience, we will always extend the solution to $[0,\infty)$ by the following definition
$$\phi(\xi):=\lim\limits_{\xi\rightarrow \ell^+}\phi(\xi),~\xi\in [0,\ell].$$

\end{prop}

Different from our use of self-similar solutions, Y. Giga, N. Mizoguchi and T. Senba obtained a blow-up profile for type I blow-up solutions of (\ref{00}) in self-similar sense.

\begin{theob}[Giga-Mizoguchi-Senba 2011, \cite{Giga-Mizo-Senba}]\label{theoa} Let $N\ge 3$. Suppose that a radial solution $(u, v )$ of (\ref{00}) undergoes type I
blowup at $t=T$ and that the blow-up set $B(u_0)=\{0\}$. Then we
have
$$(T-t)u\left((T-t)^{\frac{1}{2}}r,t\right)\rightarrow r(\Phi_\alpha)_r(r)+N\Phi_\alpha(r) \mbox{~in~} L^\infty_{loc} \mbox{~as~} t\rightarrow T,$$
where $\Phi_\alpha$ be a solution to
$$\Phi_{\xi\xi}+(\frac{N+1}{\xi}-\frac{\xi}{2})\Phi_\xi-\Phi+\Phi(\xi \Phi_\xi+N\Phi)=0 \mbox{~in~} (0,+\infty)$$
 with $\Phi^\prime(0)=0$ and $\Phi(0)=\alpha\ge 0$.
Moreover, if $(u_0)_r(r)\le 0$ is valid, we must have $\alpha>0$.

\end{theob}
\begin{rem}
Here, based on the link between blow-up profiles and self-similar solutions, we only focus on how $\phi$ decays to $0$ at infinity. Note that $\phi$ may blow up or touch zero at some $l(\neq0)$, but the estimate on its decay rate is already sufficient for our application. This leads to a significant difference between the function $\phi$ we constructed in Proposition \ref{dense S} and the function $\Phi$ in {\bf Theorem B}, although they satisfy the same equation.
\end{rem}
Another important ingredient of our proof is the zero number argument, which was initially developed by Louville and then applied to asymptotic behaviors of parabolic equations by Matano. Since the coefficients of equation (\ref{equ_w}) are not uniformly bounded in $\mathbb{R}^N$, the classical zero number argument is not applicable to our problem. Thus, we have to establish a new version for (\ref{equ_w}).
\begin{prop}\label{zero nt}
Assume that $u(r,t),v(r,t)$ are two classical solution of  (\ref{equ_w}) for $(r,t)\in [0,\infty)\times[t_1,t_2)$, and $u,v,rv_r$ are bounded. Let $U(r,t)$$=u(r,t)-v(r,t)\not\equiv0$ and $z[a,b]$ denote the zero number
in $(a, b)$. We assume in addition that there exists $L\in (0,\infty)$ such that $z_{(L,\infty]}(U(\cdot, t_1))=0$.  Then
\begin{itemize}
\item [(i)] $z_{[0,\infty)}(U(\cdot, t))<\infty$ for all $t\in (t_1, t_2)$,\\
\item [(ii)] the function $t \rightarrow z_{[0,\infty)}(U(\cdot, t))$ is nonincreasing,\\
\item [(iii)] if $U(r_0, t_0) = U_r(r_0, t_0)=0$ for some $r_0 \in [0, R]$ and $t_0 \in (t_1, t_2)$, then
$z(U(\cdot, t))<z(U(\cdot, s))$ for all $t_1<t<t_0<s<t_2.$
\end{itemize}
\end{prop}


This paper is organized as follows. In Section 2, we give out some prelininaries. In particular, we establish the zero number theorem in unbounded domain with unbounded coefficients. The main result on the blow-up profile will be shown in Section 3

\section{Preliminaries }
In this section, we introduce some notations and display some  known results. Especially, we will establish a new zero number theory which is vital for our proof.

\subsection{Local existence}

We first state the local existence and existence of blow-up solution of system (\ref{00}) which has been established in \cite{Soup-Win} (see also \cite{WinJEPE}), we repeat it here for the convenience of readers.

Rewrite problem (\ref{00}) in the form
\be{Integ eq}
    	\left\{ \begin{array}{ll}
	u_t = \Delta u-\nabla\cdot(u\nabla v),
	& x\in \mathbb{R}^N, \ t>0, \\[1mm]
	v(x,t)=\int_{\mathbb{R}^N}G(x,y)u(y,t)dy,
	& x\in \mathbb{R}^N, \ t>0,\\[1mm]
u(x,0)=u_0(x),& x\in \mathbb{R}^N,
    	\end{array} \right.
\ee
where $G$ is the Newtonian
potential. Denote by $S(t)|_{t\ge 0}$ the heat semigroup on $L^\infty (\mathbb{R}^N)$, we have the following local existence result.

\begin{prop}[Proposition 3.1 of \cite{Soup-Win}]
Let $n\ge2$ and let $u_0$ satisfy (1.3).
There exists $\tau>0$ and a unique, classical solution $(u, v)$ of (\ref{Integ eq}) such that

\be{sol space}
    	\left\{ \begin{array}{ll}
	(u, v)\in {\rm BC}^{2,1}(\mathbb{R}^N\times (0, \tau)) \times {\rm BC}^{2,0}(\mathbb{R}^N\times (0, \tau)), \\[1mm]
	u-S(t)u_0 \in {\rm BC}(\mathbb{R}^N\times (0, \tau)).
    	\end{array} \right.
\ee
Moreover, $(u, v)$ can be extended to a unique maximal solution, whose existence time $T = T(u_0, v_0)\in (0,\infty]$ satisfies
either $T>\infty$ or $\lim_{t\rightarrow T}\|u(\cdot,t)\|_{L^\infty} = \infty.$

The couple $(u, v)$ also solves (\ref{00}) and, for each $t\in (0, T)$, the function $u(\cdot, t)$ is nonnegative and radially symmetric nonincreasing. Furthermore, if we assume in addition $u_0\in L^1(\mathbb{R}^N)$, then u enjoys the mass conservation
property
$$\|u(t)\|_{L^1(\mathbb{R}^N)}=\|u_0\|_{L^1(\mathbb{R}^N)},~0<t<T.$$

\end{prop}

\subsection{Zero number theory}

In this section, we are mainly committed to establishing a new zero number theorem, that is, Proposition \ref{zero nt}. We first introduce the following zero number theorem for bounded  interval.

The zero number of a function $\psi\in C((0,R))$ is defined as the number of sign
changes of $\psi$ in $(0, R)$:
$$z_{[0,R]}(\psi) = \sup\{k\in N : {\rm ~there~are~} 0< x_0<x_1<\cdots<x_k<R
{\rm ~s.~t.~} \psi(x_i)\psi(x_{i+1})<0 {\rm ~for~} 0\le i<k\}. $$

Consider
\be{appendix}
    	\left\{ \begin{array}{ll}
	U_t = U_{rr}+\frac{n-1}{r}U_r +QU,
	& r\in (0,R), \ t\in (t_1,t_2), \\[1mm]
	U_r(0,t)=0,
	& t\in (t_1,t_2),
    	\end{array} \right.
\ee
where $U(r,t)\in C([0,R]\times [t_1,t_2])\cap W_\infty^{2,1}((0,R)\times(t_1,t_2))$ and $Q(r,t)\in L^\infty((0,R)\times(t_1,t_2))$.

\begin{lem}[Theorem 52.28 of \cite{Quit-Soup}]\label{bounded zero theory}
Assume $0<R<\infty$. Let $Q,~U$ be as above, $U\not\equiv 0$, and either $U(R, t)=0$ for all
$t\in [t_1, t_2]$ or $U(R, t)\neq 0$ for all $t \in [t_1, t_2]$. Let $z = z[0,R]$ denote the zero number
in $(0, R)$. Then
\begin{itemize}
\item [(i)] $z(U(\cdot, t))<\infty$ for all $t\in (t_1, t_2)$,\\
\item [(ii)] the function $t \rightarrow z(U(\cdot, t))$ is nonincreasing,\\
\item [(iii)] if $U(r_0, t_0) = U_r(r_0, t_0)=0$ for some $r_0 \in [0, R]$ and $t_0 \in (t_1, t_2)$, then
$z(U(\cdot, t))<z(U(\cdot, s))$ for all $t_1<t<t_0<s<t_2.$
\end{itemize}
\end{lem}

Let $u$ be a solution of \eqref{Integ eq} and $w(r,t)=r^{-n}\int_0^r s^{n-1}u(s,t)\,{\rm ds}$. Then $w$ solves the Cauchy problem (see \cite{Soup-Win} for more details)

\be{w eq}
w_t=w_{rr}+\frac{N+1}{r}w_{r}+(nw+r w_r)w,
\ee
with
\be{def of w_0}
w(r,0)=r^{-N}\int_0^r s^{N-1}u_0(s)\, {\rm ds}:=w_0(r).
\ee
Next, we will show that the zero number theory is also valid for problem \eqref{w eq}. In order to do this, we need the follwowing comparison principle.

\begin{lem}\label{CP with unbd}
Assume that $u(r,t),v(r,t)\in C^{2,1}\left([a,\infty)\times [0,T)\right)$   are two classical solution of  (\ref{w eq}) where $a>0$ is a constant. Assume that $u(r,0)< v(r,0),~x\in[a,\infty)$ and $u(a,t)< v(a,t),~t\in[0,T)$.
Assume also that $u,v,rv_r$ are uniformly bounded in $[a,\infty)\times [0,T)$.
Then $$u(r,t)<v(r,t),~ (r,t)\in [a,\infty)\times [0,T).$$
\end{lem}
\proof
Since $u,v,rv_r$ are uniformly bounded, there exists $M>0$ such that
$$|u|+|v|+|rv_r|<M,~(r,t)\in [a,\infty)\times[0,T).$$
Let $g(x,t)=u(x,t)-v(x,t)$, it is easy to show that $g$ satisfies
\be{W eq}
g_t=g_{rr}+\left(\frac{N+1}{r}+ru\right)g_{r}+(nu+nv+r v_r)g.
\ee
Let's further define $G(r,t):=e^{-Kt}g-\varepsilon e^t-\varepsilon^2 \ln (r+1)$ with $K=(2n+1)M$ and $\varepsilon<\frac{1}{\frac{N+1}{a}+M}$. Using \eqref{W eq}, we thus obtain that
\begin{align*}
G_t-G_{rr}-\left(\frac{N+1}{r}+ru\right)G_{r}&=(nu+nv+r v_r-K)e^{-Kt}g-\varepsilon e^t-\frac{\varepsilon^2}{(r+1)^2}+\left(\frac{N+1}{r}+ru\right)\frac{\varepsilon^2}{r+1} \\\nonumber
&\le (nu+nv+r v_r-K)G-\varepsilon+\varepsilon^2\left(\frac{N+1}{a}+M\right)\\\nonumber
&\le(nu+nv+r v_r-K)G,
\end{align*}
becasue $(nu+nv+rv_r-K)\le 0$.
Since $G(r,0)<0$, $G(a,t)<0$ and $\lim_{r\rightarrow\infty}G(r,t)=-\infty$, we can derive $G(r,t)\le 0$ immediately. For any fixed $(r,t)\in[a,\infty)\times [0,T)$, we can obtain $g(r,t)\le 0$ by letting $\varepsilon\rightarrow 0^+$. Moreover, we must have $g<0$ since otherwise $g(r,t)$ will touch $0$ at some pointe $(x_0,t_0)$ which is contradict to the strong comparison principle.  \qed

Now we can give the proof of Proposition \ref{zero nt}.

\proof[Proof of  Proposition \ref{zero nt}]
Thanks to the assumption $z_{[L,\infty]}(U(\cdot, t_1))=0$, we know that $U(r,t_1)\ge 0,~r\in[L,\infty)$ or $U(r,t_1)\le 0,~r\in[L_1,\infty)$. Without loss of generality, assume that the former is true.  Combining with the strong maximum principle and the assumption $U\not\equiv 0$, we can futher  choose $L_1>0$ such that $$U(r,t_1)\ge 0,~r\in[L_1,\infty) \mbox{~and~} U(L_1,t_1)> 0.$$

 By continuity of $U$, there exists $\tau>0$ such that $U(L_1,t)>0,~\forall t\in[t_1,t_1+\tau)$. On the one hand, it is easy to show that $U$ satisfies (\ref{W eq}) and then we can obtain by using Lemma \ref{CP with unbd} that $$U(r,t)>0,~\forall (r,t)\in [L_1,\infty)\times[t_1,t_1+\tau).$$
 Furthermore, we can apply Lemma \ref{bounded zero theory} on $[0,L_1]$ to derive that  $z_{[0,L_1]}(U(\cdot, t))<\infty$ for all $t\in (t_1, t_2)$ and the function $t \rightarrow z(U(\cdot, t))$ is nonincreasing in $[t_1,t_1+\tau)$. Combining these two results together, we know that $z_{[0,\infty)}(U(\cdot, t))=z_{[0,L_1]}(U(\cdot, t))<\infty$ for all $t\in (t_1, t_2)$  and the function $t \rightarrow z_{[0,\infty)}(U(\cdot, t))$ is nonincreasing in $[t_1,t_1+\tau)$. Define $t^*$ as follows
$$t_*:=\sup\{t\in (t_1,t_2)| z_{[0,\infty)}(U(\cdot, t))<\infty {\rm ~is~nonincreasing~in~} [t_1,t)\}.$$
Next we want to prove that $t_*=t_2$. Assume by contradictin that $t_*<t_2$ and denote that $k=\lim_{t\rightarrow t_*^-}z_{[0,\infty]}(U(\cdot,t))$. We know that there are two possibilities: $z_{[0,\infty]}(U(\cdot,t_*))\le k$ or $z_{[0,\infty]}(U(\cdot,t_*))> k$.

If $z_{[0,\infty]}(U(\cdot,t_*))\le k$, we can show that there exists $\tau_1>0$ such that the function $t \rightarrow z(U(\cdot, t))$ is nonincreasing in $[t_*,t_1+\tau)$ by repeating our previous discussion. Therefore, $t \rightarrow z_{[0,\infty)}(U(\cdot, t))$ is nonincreasing in $[t_1,t_*+\tau_1)$ which is contradict to the definition of $t_*$.

For the case $z_{[0,\infty]}(U(\cdot,t_*))> k$, we can choose $L_2>0$ such that $U(L_2,t_*)\neq 0$ (without loss of generality, we may assume that $U(L_2,t_*)>0$) and $z_{[0,L_2]}(U(\cdot,t_*))=k+1$. By continuity of $U$, there exists $\tau_2>0$ such that $U(L_2,t)>0,~\forall t\in [t_*-\tau_2,t_*]$. Thus, we can apply Lemma \ref{bounded zero theory} on $[0,L_2]$ to obtain that $$z_{[0,L_2]}(U(\cdot,t_*))\le \lim_{t\rightarrow t_*^-}z_{[0,L_2]}(U(\cdot,t))\le k.$$
That contradicts to the fact that $z_{[0,L_2]}(U(\cdot,t_*))=k+1$.

Therefore, we must have $t_*=t_2$ and conclusions (i) and (ii) are verified.  Based on (i) and (ii), the conclusion (iii) can be obtained by applying Lemma \ref{bounded zero theory} directly. \qed

\section{Asymptotic profile for $n\ge 3$}

In this section, we will establish the proof of Theorem \ref{thm1}.

\subsection{Construction of backward selfsimilar solution}

In this subsection, we are aiming to construct a family of backward selfsimilar solutions of (\ref{w eq}).
Although dealing with the same equation, we want to point out that our construction is obviously different from that in \cite{Giga-Mizo-Senba}: the solution constructed in \cite{Giga-Mizo-Senba} starts from $0$ and its existence interval may be bounded, while the solution we want to construct starts from the infinity. First, let $v(\xi,\tau)=(T-t)w(r,t)$ with $\xi=(T-t)^{-\frac{1}{2}}r$ and $\tau=-\ln (T-t)$.
It is routine to check that (also see formular (1.6) in \cite{Giga-Mizo-Senba}), $v$ satisfies

\be{v eq}
    	\left\{ \begin{array}{ll}
	v_\tau = v_{\xi\xi}+(\frac{N+1}{\xi}-\frac{\xi}{2})v_\xi-v+v(\xi v_\xi+Nv),
	& \xi\in(0,\infty), \ \tau\in (\tau_0,\infty), \\[1mm]
	v(\xi,\tau_0)=v_0(\xi)=Tw(T^{\frac{1}{2}}\xi,0),
	& \xi\in (0,\infty),
    	\end{array} \right.
\ee
where $\tau_0=-\ln T$.

The steady state problem related to (\ref{v eq}) is as follows.

\be{phi eq}
\phi_{\xi\xi}+\left(\frac{N+1}{\xi}-\frac{\xi}{2}\right)\phi_\xi-\phi+\phi(\xi \phi_\xi+N\phi)=0,~\xi\in(0,\infty).
\ee

Let $V(s):=\phi(\xi)\xi^2$ with $s=\frac{1}{\xi}$, then \eqref{phi eq} can be transformed into the following form.
\be{V eq}
V_{ss}+(5-N)V_s s^{-1}+\frac{1}{2}V_ss^{-3}-V_sVs^{-1}-2(N-2)Vs^{-2}+(N-2)V^2 s^{-2}=0, s\in(0,\infty).
\ee
We give the following local existence results about (\ref{V eq}).

\begin{theo}\label{loc exist of V}
There exists a subset $\mathbb{S}\subset\mathbb{R}_+$ which is dense in $\mathbb{R}_+$ such that for any $M\in \mathbb{S}$, the problem (\ref{V eq}) admits a local solution $V(r)$ fulfilling $V(0)=M$.
\end{theo}

In order to prove Theorem \ref{loc exist of V}, we need to consider an auxilary problem
\be{tilde V eq}
    	\left\{ \begin{array}{ll}
	\mathcal{L}\tilde{V}=0,
	& s>0, \\[1mm]
	\tilde{V}(0)=m,
    	\end{array} \right.
\ee
where the operator $\mathcal{L}$ is defined by $\mathcal{L}V:=(5-N)V_s s^2+\frac{1}{2}V_s -V_sVs^2-2(N-2)Vs+(N-2)V^2s.$
\begin{lem}\label{V loc lem}
For any $m>0$, there exists $L_m>0$ such that problem (\ref{tilde V eq}) admits a unique classical solution $\tilde{V}(s)\in C^2([0,L_m])$ with property $\tilde{V}_s(0)=0$ and $(2-m)\tilde{V}_{ss}(s)\ge (N-2)m(2-m)^2$.
\end{lem}

\proof
Obviously, $\mathcal{L}V=0$ is equivalent to $V^\prime=F(V,s)$ with $$F(V,s)=\frac{2(N-2)Vs-(N-2)V^2s}{\frac{1}{2}+(5-N)s^2-Vs^2}.$$  By standard ODE theory, for all $m>0$, there exists $\overline{L}_m>0$ small such that there exists a classical solution $\tilde{V}(s)$ satisfing (\ref{tilde V eq}) in the interval $(0,\overline{L}_m)$.
It is routine to check that
$$\tilde{V}_s(0)=\lim_{s\rightarrow 0^+} \frac{2(N-2)Vs-(N-2)V^2s}{\frac{1}{2}+(5-N)s^2-Vs^2}=0.$$ Moreover, for $s$ close to $0$, direct computation gives
\begin{align}\nonumber
\tilde{V}_{ss}(s)&=\frac{d}{ds}F(\tilde{V},s)\\\nonumber
&=(N-2) \frac{(2\tilde{V}s-\tilde{V}^2s)_s(\frac{1}{2}+(5-N)s^2-\tilde{V}s^2)-(2\tilde{V}s-\tilde{V}^2s)(\frac{1}{2}+(5-N)s^2-\tilde{V}s^2)_s}{[\frac{1}{2}+(5-N)s^2-\tilde{V}s^2]^2}\\\nonumber
&=(N-2)\frac{[\tilde{V}(2-\tilde{V})+o(1)][\frac{1}{2}+o(1)]-o(1)}{[\frac{1}{2}+(5-N)s^2-\tilde{V}s^2]^2}\\
&\rightarrow 2(N-2)m(2-m),~s\rightarrow 0^+.
\end{align}

For any $m\neq 2$, we have $(2-m)\tilde{V}_{ss}(s)\rightarrow 2(N-2)m(2-m)^2>(N-2)m(2-m)^2$ as $s\rightarrow 0^+$. Therefore, when $m\in (0,2)\cup(2,\infty)$, we can choose $L_m<\overline{L}_m$ such that, $(2-m)\tilde{V}_{ss}(s)\ge (N-2)m(2-m)^2$. Since $(2-m)\tilde{V}_{ss}(s)\ge (N-2)m(2-m)^2$ is automatically established, the proof is complete. \qed

For $m\in (0,2)$, we may assume the function $\tilde{V}$ constructed above satisfies $0<\tilde{V}(s)<2, s\in (0,L_m)$ by shrinking $L_m$ if neccesary. For any $\varepsilon\in(0,L_m)$, we consider another auxilary problem
\be{h eq}
    	\left\{ \begin{array}{ll}
	h_t = h_{ss}+s^{-3}\mathcal{L}h,
	& x\in(\varepsilon,L_m), \ t>0, \\[1mm]
	h(\varepsilon,t)=\tilde{V}(\varepsilon),~h(L_m,t)=\tilde{V}(L_m),
	& t>0, \\[1mm]
	h(s,0)=h_0(s)=\tilde{V}(s),
	& x\in(\varepsilon,L_m).
    	\end{array} \right.
\ee

\begin{lem}\label{h property}
Suppose $m<2$ and $0<\tilde{V}(s)<2, s\in (0,L_m)$, problem (\ref{h eq}) admits a global bounded classical solution $h_{\varepsilon,L_m}$. Moreover, the solution $h_{\varepsilon,L_m}$ (for simplicity, write it as $h$) has property
$$h_t\ge 0,~ h_s\ge 0, ~(s,t)\in (\varepsilon,L_m)\times(0,\infty).$$
\end{lem}

\proof

The fact $m\in(0,2)$ and $(2-m)\tilde{V}_{ss}(s)\ge (N-2)m(2-m)^2$ imply
$\tilde{V}_{ss}(s)>0,~s\in(\varepsilon,L_m)$. Combining with the property $\tilde{V}_s(0)=0$, we can obtain that
\begin{equation}\label{V s}
\tilde{V}_{s}(s)>0,~s\in(\varepsilon,L_m).
\end{equation}
 Next, we will prove the main conclusions through three steps.

{\bf Step 1. Global existence and boundedness.}

Since the problem \eqref{h eq} is regular in the interval $(\varepsilon,L_m)$ , the local existence and uniqueness of solution (with maximum existence time $T_h$) can be obtained by standard parabolic theory. By maximum principle, we know that $h$ is positive. It is easy to know that $2$ is an upper solution of problem \eqref{h eq}, so we can derive $h<2,~x\in(\varepsilon,L_m),t\in [0,\infty)$ by comparison principle.  Therefore, we know that $h$ is a global solution and it is bounded.

{\bf Step 2. Monotonicity: $h_t\ge0$.}

Direct computation gives that
 $$\tilde{V}_{ss}(s)+s^{-3}\mathcal{L}\tilde{V}\ge (N-2)m(2-m)>0,~s\in(\varepsilon,L_m).$$
By comparing $h$ with the subsolution $\underline{h}(s,t) := \tilde{V}(s)$ via comparison principle, we
obtain  $h(s,\tau)\ge \tilde{V}(s)$ in $(\varepsilon,L_m)\times[0,\infty)$.

Fix $\tau>0$ and set $g(t):= h(t+\tau)$. Since $g(s,0)=h(s,\tau)\ge h(s,0)$, we infer
from comparison principle that $g\ge h$ on $(0,\infty)$. The conclusion $h_t\ge0$ then follows by the definition of derivative.

{\bf Step 3. Monotonicity: $h_s\ge0$.}

Let $H(s,t)=h_s(s,t)$, taking derivative on both side of the equation of $h$, we can derive
\begin{align}\nonumber
H_t&=H_{ss}+[(5-N)s^{-1}+\frac{1}{2}s^{-3}-Vs^{-1}]H_s+[-(5-N)s^{-2}-\frac{3}{2}s^{-4}-Hs^{-1}+Vs^{-2}]H\\\nonumber
&+[-2(N-2)s^{-2}+2(N-2)Vs^{-2}]H+2(N-2)s^{-3}[2V-V^2]\\\label{Heq}
&\ge H_{ss}+p(s)H_s+q(s)H,
\end{align}
with $$p(s)=(5-N)s^{-1}+\frac{1}{2}s^{-3}-Vs^{-1}$$ and $$q(s)=(-1-N)s^{-2}-\frac{3}{2}s^{-4}-Hs^{-1}+(2N-3)Vs^{-2}.$$
Using the property $h_t\ge 0$ and \eqref{V s}, we can derive
$$h(s,t)\ge h(s,0)=\tilde{V}(s)\ge \tilde{V}(\varepsilon)= h(\varepsilon,t).$$
It can be obtained by the definition of derivative
\begin{equation}\label{lbdr}
H(\varepsilon,t)=h_s(\varepsilon,t)\ge 0.
\end{equation}
The fact $\tilde{V}(L_m)<2$ implies that
\begin{equation}\label{LV}
\mathcal{L}\left[\tilde{V}(L_m)\right]=-(N-2)\tilde{V}(L_m)L_m(2-\tilde{V}(L_m))<0.
\end{equation} Thanks to \eqref{V s} and \eqref{LV}, it is easy to check that $\overline{h}(s,t):\equiv \tilde{V}(L_m)<2$ is a super-solution of (\ref{h eq}). By comparison principle, we have $h(s,t)\le \tilde{V}(L_m)=h(L_m,t)$ which furthermore implies
\begin{equation}\label{rbdr}
H(L_m,t)=h_s(L_m,t)\ge 0.
\end{equation}

Combining \eqref{V s}, \eqref{Heq}, \eqref{lbdr} and \eqref{rbdr}, an comparison argument shows that
$$H(s,t)\ge 0, s\in(\varepsilon,L_m),~t>0.$$
The proof is complete.\qed

For $m>2$, we may assume the function $\tilde{V}$ constructed in Lemma \ref{V loc lem} satisfies $2<\tilde{V}(s)<m, s\in (0,L_m)$ by shrinking $L_m$ if neccesary.  A similar argument as in proof of Lemma \ref{h property} shows that
\begin{lem}\label{h property2}
Suppose $m>2$ and $2<\tilde{V}(s)<m, s\in (0,L_m)$, problem (\ref{h eq}) admits a global bounded classical solution $h_{\varepsilon,L_m}$ (write it as $h$ for simplicity). Moreover, the solution has property
$$h_t\le 0,~ h_s\le 0, ~(s,t)\in (\varepsilon,L_m)\times(0,\infty).$$
\end{lem}

Now, we turn to the proof of Theorem \ref{loc exist of V}.

\proof[Proof of Theorem \ref{loc exist of V}]

To verify the conclusion, fix $m\in (0,2)\cup(2,\infty)$, for any $\delta>0$ we construct a function $V\in C([0,L_m])\cup C^2(0,L_m)$ fulfilling problem (\ref{v eq}) with initial data $V(0)\in [m-\delta,m+\delta]$. We will only deal with the case $m\in(0,2)$, the other case $m\in(2,\infty)$ can be treated similarly. For any fixed $m\in (0,2)$ and arbitrary $\delta>0$ small, we can shrink $L_m$ such that $\tilde{V}(L_m)<m+\delta.$

For any $a\in (0,L_m)$, we know that $h_{\varepsilon,L_m}(a,t)$ is well-defined for all $\varepsilon\in(0,a)$, we show that $h_{\varepsilon,L_m}(a,t)$ is decreasing with respect to $\varepsilon$ in $(0,a)$. Assume that $0<\varepsilon_1<\varepsilon_2<a$. Thanks to Lemma \ref{h property}, we know that  $h_{\varepsilon_i,L_m},i=1,2$ is nondecreasing with respect to both of its variable and $\tilde{V}$ is also nondecreasing, which means that
$$h_{\varepsilon_1,L_m}(\varepsilon_2,t)\ge h_{\varepsilon_1,L_m}(\varepsilon_2,0)=\tilde{V}(\varepsilon_2)=h_{\varepsilon_2,L_m}(\varepsilon_2,t).$$

In view of $h_{\varepsilon_1,L_m}(L_m,t)=\tilde{V}(L_m)=h_{\varepsilon_1,L_m}(L_m,t)$, a comparison argument shows that
 $$h_{\varepsilon_1,L_m}(s,t)\ge h_{\varepsilon_2,L_m}(s,t),~s\in [\varepsilon_2,L_m],~t>0.$$
Choosing $s=a$ yields the conclusion.

Therefore, by the properties $h_t\ge 0, h\leq \tilde{V}(L_m)\leq 2$ and the interior regularity of $h$, via a diagonal argument we can deduce
$h_{\varepsilon_1,L_m}(s,t)\rightarrow V(s)$ pointwise in $(0,L_m]$ and uniformly in any compact subset of $(0,L_m]$ as $\varepsilon\rightarrow 0^+$. Moreover, $V(s)$ satisfies equation (\ref{v eq}) in $(0,L_m)$.

Since $h(s,t)$ is increasing in $s$, so as $V(s)$. Thus, we have $V(s)\leq V(L_m)=\lim_{t\rightarrow\infty}h(L_m,t)=\tilde{V}(L_m)<m+\varepsilon.$ Moreover, we also have $m=h(s,0)\le h(s,t)\le V(s)$. Thus, we establish the following estimate $m\le V(s)\le m+\delta$. Since $V(s)$ is increasing in $s$, we can define $V(0):=\lim_{s\rightarrow 0^+}V(s)\in [m,m+\delta]$. The proof is complete. \qed

Turn back to $\phi$, for $m\in \mathcal{S}$, we get a function $\phi(\xi), \xi\in (l,\infty)$ with $l=\frac{1}{L_m}$, satisfying

\be{phi eq2}
\phi_{\xi\xi}+(\frac{N+1}{\xi}-\frac{\xi}{2})\phi_\xi-\phi+\phi(\xi \phi_\xi+N\phi)=0,
\ee
and
$$\phi(\xi)\xi^2\rightarrow m,~\xi\rightarrow\infty.$$

Extend $l>0$ as small as possible to a minimal value, one says
\be{L*}
\ell:=\inf\{l>0|~ \phi {\rm ~can~be~defined~in~} (l,\infty)\}.
\ee

We have the following result.

\begin{lem}\label{prop phi}
Assume that $\ell$ is defined as in (\ref{L*}) and $\phi(r)$ is the corresponding positive solution of (\ref{phi eq2}). Then precisely one of the three alternatives holds:
\begin{itemize}
\item $\phi$ is unbounded: $\phi(\xi)\rightarrow\infty$ as $\xi\rightarrow \ell^+$; \\
\item $\phi$ is bounded and $\ell>0$: $\lim\limits_{\xi\rightarrow \ell^+}\phi(r)=0$;\\
\item $\phi$ is bounded and $\ell=0$: $\lim\limits_{\xi\rightarrow \ell^+}\phi(\xi)$ exists and $\phi^\prime(\xi)\rightarrow 0$ as $\xi\rightarrow \ell^+$.
\end{itemize}
\end{lem}

\proof

{\bf Case 1. $\phi$ is unbounded near $\ell$.}

First, we show that $\phi$ is nonincreasing in some right neighborhood of $\ell$ which implies $\phi(\xi)\rightarrow\infty$ as $\xi\rightarrow \ell$ immediately. If the decreasing property is not true, then we have the following claim.

{\bf Claim:} There exists a decreasing sequence $\{r_n\}_{n=1}^\infty$ that converges to $\ell$ such that $\phi(r_n)>\frac{1}{N}$ and $\phi^\prime(r_n)>0$.

 Let's verify the claim. Since $\phi$ is unbounded, there exists a decreasing sequence $\{s_n\}_{n=1}^\infty$ such that $s_n\rightarrow \ell$ and $\phi(s_n)>\frac{1}{N}$. For any $n$, we determine $r_n$ as follows. If $\phi^\prime(s_n)>0$, let $r_n:=s_n$. If $\phi^\prime(s_n)\le 0$, consider
$$s_n^*=\sup\{s<s_n|\phi^\prime(s)>0\}.$$  It is mentioned that $\{s<s_n|\phi^\prime(s)>0\}=\emptyset$ can not be ture since it will leads to the boundedness of $\phi$. It is routine to check that $\phi(s_n^*)>\frac{1}{N}$ and there exists $\delta\in (0,s_n^*)$ such that $\phi(s)>\frac{1}{N}, s\in (s_n^*-\delta,s_n^*)$. Thanks to the definition of $s_n^*$ we have $(s_n^*-\delta,s_n^*)\cup \{r|\phi^\prime(r)>0\}\neq\emptyset$. Therefore we can choose arbitrary $r_n\in (s_n^*-\delta,s_n^*)\cup \{r|\phi^\prime(r)>0\}$.

Next, we will show that this claim leads to a contradiction. Let $\{r_n\}_{n=1}^\infty$ be as in the claim, we define
$$\overline{r}_n:=\sup\{r\in(r_n,\infty)|~\phi^\prime(s)>0,~s\in (r_n,r)\}.$$
If $\overline{r}_n\leq \sqrt{N+1}$ is valid for some $n$, we have $\phi^\prime(r)>0,~r\in (r_n,\overline{r}_n)$ and $\phi^\prime(\overline{r}_n)=0$. Moreover, it is routine to check that $\phi(r)>\phi(r_n)\ge \frac{1}{N},~r\in (r_n,\overline{r}_n)$ and $\frac{N+1}{r}-\frac{r}{2}+r\phi\ge \frac{N+1}{2r},~r\in (r_n,\overline{r}_n)$.
Therefore, for $r\in (r_n,\overline{r}_n)$ we have
$$\phi^{\prime\prime}+\frac{N+1}{2r}\phi^\prime\le  \phi^{\prime\prime}+\left(\frac{N+1}{r}-\frac{r}{2}+r\phi\right)\phi^\prime= \phi(1-N\phi)\leq 0.$$
As a result, we obtain $$\left(r^{\frac{N+1}{2}}\phi^\prime\right)^\prime\leq 0.$$
Integrating from $r_n$ to $\overline{r}_n$ we can further derive
$$0<r_n^{\frac{N+1}{2}}\phi^\prime(r_n)\leq r_n^{\frac{N+1}{2}}\phi^\prime(\overline{r}_n)=0,$$
which is absurd.

Thus, we must have $\overline{r}_n>\sqrt{N+1}$ for all $n\in \mathbb{N}$. Hence, we must have $\phi^\prime(r)>0,~r\in(r_n,\sqrt{N+1})$. Letting $n\rightarrow\infty,$ we get $\phi^\prime(r)>0,~r\in(\ell,\sqrt{N+1})$ which is contradict to the unboundedness of $\phi$ (near $\ell$).

{\bf Case 2. $\phi$ is bounded and $\ell>0$.}

Define $a(r)=e^{-\frac{r^2}{4}+\int_0^r \rho \phi(\rho)\ {\rm d\rho}}$. Direct computation gives

\begin{align}\label{r a}
\nonumber\left(r^{N+1}a(r)\phi^\prime(r)\right)^\prime&=r^{N+1}a(r)\left[\phi^{\prime\prime}+\left(\frac{N+1}{r}-\frac{r}{2}+r\phi\right)\phi^\prime\right]\\
&= r^{N+1}a(r)\phi(1-N\phi).
\end{align}

It easy to show that $\phi$ is uniformly continuous near $\ell$ by the boundedness $\phi$ and \eqref{r a}, which implies that $\phi(\ell):=\lim_{r\rightarrow\ell^+}\phi(r)\ge 0$ is well-defined. If $\phi(\ell)> 0$, we can extend the solution a little more to $(\ell-\delta_1,\infty)$ ($\delta_1>0$ is some small constant) which is contradict to the definition of $\ell$.

{\bf Case 3. $\phi$ is bounded and $\ell=0$.}

Let $s_n=\frac{1}{n}$, from the boundedness of $\phi$, it is obvious that $s_n^{N+1}\phi(s_n)\rightarrow 0$.

Integrating (\ref{r a}) from $s_n$ to $\xi$ yields
$$\xi^{N+1}a(\xi)\phi^\prime(r)-s_n^{N+1}a(s_n)\phi^\prime(s_n)=\int_{s_n}^\xi \rho^{N+1}a(\rho)\phi(1-N\phi)\ {\rm d\rho}.$$
Letting $n\rightarrow \infty$, we have
$$\phi^\prime(\xi)=\frac{1}{a(\xi)}\int_0^\xi\left(\frac{\rho}{\xi}\right)^{N+1}a(\rho)\phi(1-N\phi)\ {\rm d\rho},$$
which implies $\phi^\prime(\xi)\rightarrow 0^+$ as $\xi\rightarrow 0^+$.

The conclusion $\phi^\prime(\xi)\rightarrow 0,~\xi\rightarrow 0^+$ ensures that $\phi$ is uniformly continuous near $0$, thus the limit $\lim\limits_{\xi\rightarrow 0^+}\phi$ must exist.\qed

\subsection{Asymptotic profile of $W$}

In this subsection, we are committed to establishing the exact blow-up profile of $W$. The main result can be stated as follows.

\begin{theo}\label{thm3.1}
Suppose that $u_0$ satisfies (\ref{initial 1}) and (\ref{initial w_0}),  $w$ is a single point blow-up solution of (\ref{w eq}) with $T<\infty$. Then the following limit exists
$$\lim_{r\rightarrow 0}W(r)r^2\in [0,\infty].$$
\end{theo}

Let $\phi_m$ be the solution constructed in last subsetion with property $\lim_{\xi\rightarrow \infty}\phi_m(\xi)\xi^2=m$, we will give the proof of Theorem \ref{thm3.1} via a zero number discussion.


\proof[Proof of Theorem \ref{thm3.1}.]

Assume by contradiction that
$$\limsup_{r\rightarrow0^+} W(r)r^2>\liminf_{r\rightarrow0^+} W(r)r^2.$$
Therefore, we can choose a decreasing sequence $\{r_n\}_{n=1}^\infty$ such that $$\lim_{n\rightarrow\infty}W(r_{2n})r_{2n}^2=A>B=\lim_{n\rightarrow\infty}W(r_{2n+1})r_{2n+1}^2.$$
Since $w_0$ satisfies \eqref{initial w_0} and $\mathcal{S}$ is dense in  $\mathbb{R}^+$, we can choose $m\in \mathcal{S}\cap (\frac{3A+B}{4},\frac{3B+A}{4})$ such that $\lim_{r\rightarrow\infty}w_0(r)r^2\neq m$. Let $\phi_m(\xi),~\xi\in[\ell,\infty)$ be the related solution constructed in Proposition \ref{dense S} and set $U(r,t)=w(r,t)-\frac{1}{T-t}\phi(r/\sqrt{T-t})$, we want to apply zero number theory to $U$ to get a contradiction.

The assumption $\lim_{r\rightarrow\infty}w_0(r)r^2\neq m$ insures that we can choose $L$ larege such that $z_{(L,\infty]}(U(\cdot, 0))=0$ and $U(L,0)\neq 0$ which implies that $U(r,\tau)$ has at most finite zero number point for $\tau>0$ small. We denote the maximum of the zero number of $U(r,\tau)$ by $k$.  For $t<T$, $w(r,t)$ is positive, bounded and satisfies $w_r(0,t)=0$. Therefore, we can apply Proposition \ref{zero nt} to obtain that $w(r,t)-\frac{1}{T-t}\phi(r/\sqrt{T-t}),~0<t<T$ has at most $k$ zero number in $(0,\infty)$.

Next, choose $N$ large enough such that $$W(r_{2n})r^2_{2n}>\frac{3A+B}{4},  W(r_{2n+1})r^2_{2n+1}<\frac{3B+A}{4},~\forall n\ge N.$$ Since $w(r,t)\rightarrow W(r)$ as $t\rightarrow T$, we can go a step further to choose $T_1<T$ such that
 \be{w osi}w(r_{2n},t)r^2_{2n}>\frac{3A+B}{4}, w(r_{2n+1},t)r^2_{2n+1}<\frac{3B+A}{4},~\forall t\ge T_1, n=N,N+1,\cdots,N+k.
\ee
For $r\in [r_{2N+2k+1},r_{2N}]$, we have $\frac{r}{\sqrt{T-t}}\rightarrow\infty$ uniformly as $t\rightarrow T$. Thus, direct computation shows that $$\lim_{t\rightarrow T}\frac{1}{(T-t)}\phi(r/\sqrt{T-t})r^2=\lim_{\xi\rightarrow\infty}\phi(\xi)\xi^2 =m\in  (\frac{3A+B}{4},\frac{3B+A}{4})$$ uniformly in $ [r_{2N+2k+1},r_{2N}]$.

Therefore, there exists $T_2<T$ such that
\be{interval}
\frac{r_{2N+2k+1}}{\sqrt{T-T_2}}>L,
\ee
and
\be{phi stable}
\frac{3A+B}{4}\le \frac{1}{(T-t)}\phi(r/\sqrt{T-t})r^2 \le\frac{3B+A}{4}, t>T_2, r\in [r_{2N+2k+1},r_{2N}].
\ee

Combining (\ref{w osi}), (\ref{interval}) and (\ref{phi stable}),  it is easy to conclude that for any $t>\max\{T_1,T_2\}$, $U(r,t)=w(r,t)-\frac{1}{T-t}\phi_m(r/\sqrt{T-t})$ has at least $2k$ zero number in $[L\sqrt{T-t},\infty)$, which is a contradiction. \qed

\subsection{Asymptotic profile of $U$}

In this subsection, we will verify Theorem  \ref{thm1} by establishing the exact blow-up profile of $U$. Specifically, we will prove the following results.

\begin{theo}\label{thm4.1} Under the assumptions of Theorem \ref{thm1}, we have
 $$\lim_{r\rightarrow0}U(r)r^2=\alpha,$$
where $\alpha:=(N-2)\lim_{r\rightarrow0}W(r)r^2<\infty$.
\end{theo}

\proof

Similar to the proof of Theorem  in \cite{Soup-Win}, we can obtain $w(r,t)\le \frac{C}{r^2}$ for some $C\in (0,\infty)$. Thanks to this estimate, we may employ Theorem \ref{thm3.1} to see that $\lim_{r\rightarrow0}W(r)r^2$ is finite and so does $\alpha$.

Next, we turn to the proof of $\lim_{r\rightarrow0}U(r)r^2=\alpha.$

{\bf Step 1.}  $w_{rr}\geq -\frac{C}{r^4}$.

For any fixed $t^*\in (0,T)$ and $\varepsilon>0$ small,
we choose
$$ \zeta(r)=\left\{
\begin{array}{lll}
0, & r\in [0,\frac{\varepsilon}{3}]\cup [4\varepsilon, \infty),\\
{\rm smooth~connection}, &  r\in (\frac{\varepsilon}{3},\frac{\varepsilon}{2})\cup (3\varepsilon, 4\varepsilon),\\
1, & r\in [\frac{\varepsilon}{2},3\varepsilon],
\end{array}
\right.
$$
and
$$ \eta(t)=\left\{
\begin{array}{lll}
0, & t\in [0,t^*-2\varepsilon^2],\\
{\rm smooth~connection}, &  t\in (t^*-2\varepsilon^2,t^*-\varepsilon^2),\\
1, & r\in [t^*-\varepsilon^2,\infty),
\end{array}
\right.
$$
such that $0\leq \zeta,\eta\leq 1$ and $\varepsilon^2|\eta^\prime|+\varepsilon^2|\zeta^{\prime\prime}|+\varepsilon|\zeta^\prime|\leq C.$

Let $v(s,\tau)=\eta(t)\zeta(r)w(r,t)$ with $s=\frac{r}{\varepsilon},~\tau=\frac{t-t^*}{\varepsilon^2}$.
Direct computation gives
\be{vs eq}
v_\tau=v_{ss}+f
\ee
where $f(s,\tau)=\varepsilon^2\{\eta^\prime\zeta w-2\eta \zeta^\prime w_r-\eta\zeta^{\prime\prime}w+\eta\zeta[\frac{N+1}{r}w_r+(nw+rw_r
)w]\}.$

Recall that $r^Nw(r,t)=\int_0^r s^{N-1}u(s,t)ds$.
Take derivative on both side, we obtain $r^{N-1}u(r,t)=r^Nw_r(r,t)+Nr^{N-1}w(r,t)$ which is equivalent to
\be{u of w}
u(r,t)=rw_r(r,t)+Nw(r,t).
\ee
Combine with $0<u\le \frac{C}{r^2}$, we get
\be{w wr est}
w=r^{-N}\int_0^r s^{N-1}u(s,t)ds\leq \frac{C}{(N-2)r^2},
\ee
and
$$|w_r|\leq \left|\frac{u}{r}-\frac{Nw}{r}\right|\le\frac{3C}{r^3}.$$

It is routine to check that $\|f\|_{L^\infty((\frac{1}{3},4)\times (-2,2))}\leq \frac{C}{\varepsilon^2}$ which implies $\|f\|_{L^p((\frac{1}{3},4)\times (-2,2))}\leq \frac{C}{\varepsilon^2}$ for any $p\in (1,\infty)$.
Applying standard parabolic estimate to (\ref{v eq}), we can derive
$$\|v\|_{W^{2,1}_p\left((\frac{1}{3},4)\times (-2,2)\right)}\leq C\left(\|v\|_{L^p\left((\frac{1}{3},4)\times (-2,2)\right)}+\|f\|_{L^p\left((\frac{1}{3},4)\times (-2,2)\right)}\right)\leq \frac{C}{\varepsilon^2}.$$
For $p$ large enough, by Sobolev imbedding Theorem, there exists $\alpha\in (0,1)$ such that
\be{v 1+alpha}
\|v\|_{C^{1+\alpha,\frac{1+\alpha}{2}}\left((\frac{1}{3},4)\times (-2,2)\right)}\leq \frac{C}{\varepsilon^2}.
\ee
Especially, we obtain

\be{w 1+alpha}
\|w\|_{C^{1+\alpha,\frac{1+\alpha}{2}}\left((\frac{\varepsilon}{2},3\varepsilon)\times (t^*-\varepsilon^2,t^*+\varepsilon^2)\right)}\leq \frac{C}{\varepsilon^{3+\alpha}}.
\ee

Choose $$\tilde{ \zeta}(r)=\left\{
\begin{array}{lll}
0, & r\in [0,\frac{\varepsilon}{2}]\cup [3\varepsilon, \infty)\\
{\rm smooth~connection}, &  r\in (\frac{\varepsilon}{2},\varepsilon)\cup (2\varepsilon, 3\varepsilon) \\
1, & r\in [\varepsilon,2\varepsilon]
\end{array}
\right.
$$
and
$$ \tilde{\eta}(t)=\left\{
\begin{array}{lll}
0, & t\in [0,t^*-\varepsilon^2]\\
{\rm smooth~connection}, &  t\in (t^*-\varepsilon^2,t^*) \\
1, & r\in [t^*,\infty)
\end{array}
\right.
$$
such that $0\leq \tilde{ \zeta},\tilde{ \eta}\leq 1$ and $\varepsilon^2|\tilde{ \eta}^\prime|+\varepsilon^2|\tilde{ \zeta}^{\prime\prime}|+\varepsilon|\tilde{ \zeta}^\prime|\leq C.$

Let $\tilde{v}(s,\tau)=\tilde{\eta}(t)\tilde{\zeta}(r)w(r,t)$ with $s=\frac{r}{\varepsilon},~\tau=\frac{t-t^*}{\varepsilon^2}$.
Similar as before, we have
\be{tilde v eq}
\tilde{v}_\tau=\tilde{v}_{ss}+\tilde{f}
\ee
where $\tilde{f}(s,\tau)=\varepsilon^2\left\{\tilde{\eta}^\prime\tilde{\zeta} w-2\tilde{\eta} \tilde{\zeta}^\prime w_r-\tilde{\eta}\tilde{\zeta}^{\prime\prime}w+\tilde{\eta}\tilde{\zeta}\left[\frac{N+1}{r}w_r+(nw+rw_r
)w\right]\right\}.$

Thanks to (\ref{w 1+alpha}), we can obtain $$\|\tilde{f}\|_{C^{\alpha,\frac{\alpha}{2}}\left((\frac{1}{2},3)\times(-1,1)\right)}\leq \frac{C}{\varepsilon^{2}}.$$

By Schauder eastimate, we can derive

$$\|\tilde{v}\|_{C^{2+\alpha,1+\frac{\alpha}{2}}\left((\frac{1}{2},3)\times(-1,1)\right)}\leq \frac{C}{\varepsilon^{2}}.$$

Especially, we have
$$|\partial_{ss} w(\varepsilon s,\varepsilon^2\tau)|\leq \frac{C}{\varepsilon^{2}},~\forall (s,\tau)\in [0,1]\times(1,2).$$

Equivalently, we have
$$|\partial_{rr} w(r,t)|\leq \frac{C}{\varepsilon^{4}},~\forall (r,t)\in [t^*,t^*+\varepsilon^2]\times(\varepsilon,2\varepsilon).$$

Therefore, we get
$$|w_{rr}(r,t^*)|\le\frac{C}{\varepsilon^4}, r\in [\varepsilon,2\varepsilon].$$

Since $t^*$ and $\varepsilon$ are arbitrary, we get our conclusion.

{\bf Step 2.} $u_r(r,t)\ge -\frac{C}{r^3}.$

Take derivative on both side of (\ref{u of w}) and combine with the estimate in Step 1, we get
$$u_r=rw_{rr}+(N+1)w_r\ge -\frac{C}{r^3}.$$

{\bf Step 3.} $U(r)r^2\rightarrow\alpha$ as $r\rightarrow 0$.

Assume by contradiction that there exists a decreasing sequence $\{r_n\}_{n=1}^\infty$ such that $r_n\rightarrow 0^+$ as $n\rightarrow\infty$ and $U(r_n)r_n^2\rightarrow \alpha^*\neq \alpha$.

We claim that there exists $C>0$ independent of $t$ and $n$ such that
\be{u-u}
r_n\int_{\sigma r_n}^{r_n} \left(\frac{\rho}{r_n}\right)^{N-1}\left[u(\rho,t)-u(r_n,t)]\right] d\rho\le C(1-\sigma)^2,~\sigma\in(0,2).
\ee

In fact, for $\sigma\le1$, we can calculate by the estimate of Step 2 that
\begin{align*}
{\rm LHS}&=
r_n\int_{\sigma r_n}^{r_n} \left(\frac{\rho}{r_n}\right)^{N-1}\left[\int_{\rho}^{r_n}-u_r(r,t)dr\right] d\rho\\
&\le r_n\int_{\sigma r_n}^{r_n} \left(\frac{\rho}{r_n}\right)^2\left[\int_{\rho}^{r_n}\frac{C}{r^3}dr\right] d\rho\\
&= r_n\int_{\sigma r_n}^{r_n} \left(\frac{\rho}{r_n}\right)^2\left[\frac{C}{2\rho^2}-\frac{C}{2r_n^2}\right] d\rho\\
&=\frac{C}{2}\int_\sigma^1 (1-\tilde{\rho}^2)d\tilde{\rho}\\
&\leq C(1-\sigma)^2.
\end{align*}

If $\sigma>1$, we can make similar estimates as follows
\begin{align*}
{\rm LHS}&=
-r_n\int_{ r_n}^{\sigma r_n} \left(\frac{\rho}{r_n}\right)^{N-1}\left[\int_{r_n}^{\rho}u_r(r,t)dr\right] d\rho\\
&\le r_n\int_{ r_n}^{\sigma r_n} \left(\frac{\rho}{r_n}\right)^{N-1}\left[\int_{r_n}^{\rho}\frac{C}{r^3}dr\right] d\rho\\
&= r_n\int_{ r_n}^{\sigma r_n} \sigma^{N-3}\left(\frac{\rho}{r_n}\right)^2\left[\frac{C}{2r_n^2}-\frac{C}{2\rho^2}\right] d\rho\\
&=\frac{C\sigma^{N-3}}{2}\int^\sigma_1 (\tilde{\rho}^2-1)d\tilde{\rho}\\
&\leq C(1-\sigma)^2.
\end{align*}
Letting $t\rightarrow T^-$ and after straighforward rearrangement leads to
\be{upper est}
r_n\int_{\sigma r_n}^{r_n} \left(\frac{\rho}{r_n}\right)^{N-1}U(\rho) d\rho \le r_n\int_{\sigma r_n}^{r_n} \left(\frac{\rho}{r_n}\right)^{N-1}U(r_n) d\rho+C(1-\sigma)^2.
\ee
By the definition of $W$, we have the following identity
$$\int_{\sigma r_n}^{r_n} \rho^{N-1}U(\rho)d\rho=r_n^NW(r_n)-(\sigma r_n)^N W(\sigma r_n).$$
Therefore, (\ref{upper est}) is equivalent to
$$r_n^2 W(r_n)-\sigma^{N-2}(\sigma r_n)^2 W(\sigma r_n)\leq \int_\sigma^1 s^{N-1}ds r_n^2 U(r_n)+C(1-\sigma)^2.$$
Taking $n\rightarrow\infty$, we can see that
$$(1-\sigma^{N-2})\frac{\alpha}{N-2}\le \alpha^* \int_\sigma^1 s^{N-1}ds+C(1-\sigma)^2,$$
becasue $\lim_{r\rightarrow0}W(r)r^2=\frac{\alpha}{N-2}$
due to the definition of $\alpha$.

For $\sigma<1$, divide both side by $1-\sigma$ and then let $\sigma\rightarrow 1^-$ we can derive
$$\alpha\le \alpha^*.$$

Choose $\sigma>1$, we can obtain $\alpha\ge \alpha^*$ similarly and then we get a contradiction that $\alpha=\alpha^*$.\qed

\begin{rem}
The following equivalence relation is always true: $$\lim_{r\rightarrow0}U(r)r^2=\alpha\Longleftrightarrow\lim_{r\rightarrow0}W(r)r^2=\frac{\alpha}{N-2}.$$
In fact, one side is just proved in Theorem \ref{thm4.1} and the other side is trivial. Moreover, based on the process of above proof, this result is also valid for bounded domain(radially symmetric case).

It can be seen from the proof process of \cite{Soup-Win} and the proof of the our main result, the assumption ``$u_0$ is nonincreasing in $r=|x|$" in \eqref{initial 1} can be weaken into ``$w_0=r^{-N}\int_0^r s^{N-1}u_0(s)ds$ is nonincreasing in $r$". Obviously, $u_0$ is nonincreasing in $r$ implies $w_0$ is nonincreasing in $r$. But the opposite is not true. To see this, let $\phi_2(r),~r\in [2,\infty)$ be any nonnegative nonincreasing smooth function, we choose $\phi_1(r),~r\in [1,2]$ satisfies
$$\phi_1(1)=\phi_1^\prime(1)=\phi_1^{\prime\prime}(1)=(\phi_1-\phi_2)(2)=(\phi_1-\phi_2)^\prime(2)=(\phi_1-\phi_2)^{\prime\prime}(2)=0,$$
and construct
$$ u_0(r)=\left\{
\begin{array}{lll}
kr^2(r-1)^4, & r\in [0,1],\\
\phi_1(r), &  r\in (1,2), \\
\phi_2(r), & r\in [2,\infty).
\end{array}
\right.
$$
It is routine to check that $w_0=r^{-N}\int_0^r \rho^{N-1}u_0(\rho)\, {\rm d\rho}$ is nonincreasing in $r$ for large $k$, but $u_0$ is not nonincreasing in $r$.
\end{rem}

\end{document}